\documentclass[aps,pre,preprint]{revtex4}
\usepackage{amsmath,bm,epsfig}

\usepackage{amsmath}
\usepackage{amsfonts}
\usepackage{amssymb}
\usepackage{upgreek}

\def\Integer{\mathbb{Z}}

\def\be{\begin{equation}}\def\ee{\end{equation}}
\def\bea{\begin{eqnarray}}\def\eea{\end{eqnarray}}

\usepackage{graphicx}
\usepackage{amsmath,amsfonts,amsthm}
\usepackage{amssymb}
\usepackage{comment}
\usepackage{enumitem}
\usepackage[utf8]{inputenc}
\usepackage[T1]{fontenc}
\usepackage{emptypage}
\usepackage{hyperref}
\usepackage{color}
\usepackage{mathrsfs}

\newcommand{\Z}{\mathbb{Z}}
\newcommand{\R}{\mathbb{R}}
\newcommand{\ki}{\mathbf{k}}

\begin{document}

\title{Explicit form of the effective evolution equation for  \\the randomly forced Schr\"{o}dinger equation with quadratic nonlinearity}

\author{Huilin Zhang$^{1,2}$, Elena Tobisch$^2$}

\email{Elena.Tobisch@jku.at}

\affiliation{$^1$ School of Mathematics, Fudan University, Shanghai, China\\
$^2$ Institute of Analysis, Johannes Kepler University, Linz, Austria\\
}
% $^3$University Paris 7, Paris, France}

%{PACS: 05.45.-a, 43.25.+y, 47.35.Jk}
\begin{abstract}
 An effective equation describes a weakly nonlinear wave field evolution governed by  nonlinear dispersive PDEs \emph{via} the set of its resonances in  an arbitrary big but finite domain in the Fourier space. We consider the   Schr\"{o}dinger equation with quadratic nonlinearity including small external random forcing/dissipation.
  An effective equation is deduced explicitly for each  case of monomial quadratic nonlinearities $ u^2, \, \bar{u}u, \, \bar{u}^2$ and  the sets of  resonance clusters are studied. In particular, we demonstrate that the  nonlinearity $\bar{u}^2$ generates no 3-wave resonances and its effective equation is degenerate while in two other cases the sets of resonances are not empty.
Possible  implications for wave turbulence theory are briefly discussed.

\end{abstract}

\maketitle
%\tableofcontents
\section{Introduction}
The study of a nonlinear wave field evolution under the action of a  small random external force is met in many physical applications including, but not restricted to, oceanology, nonlinear optics, physics of atmosphere. The traditional approach for studying statistical properties of this type in weakly nonlinear wave systems is called  kinetic wave turbulence theory (WTT) and  originates from \cite{zak2}, while the current state of the art can be found in \cite{NR11}. Under  a set of additional assumptions, a wave kinetic equation is deduced  and the study of its stationary solutions is the main subject of the theory. If one of the fundamental assumptions  of the theory - infiniteness of physical space - is omitted, the classical wave kinetic equation becomes invalid. Instead, the dynamics of resonance clusters, if they exist, describes the total nonlinear field evolution at suitable time scale, as it was first demonstrated in \cite{PRL}, and this is the subject of discrete WTT.

Quite recently, a novel approach has been progressed  which allows to develop a unifying description of both discrete and statistical properties of these wave systems, e.g. \cite{YaYo2013,KM2015,FGH2016}. It  describes limiting behavior of small amplitude solutions for a nonlinear dispersive PDE in a finite volume by preserving only its resonant terms, the resulting equation is called below  the effective equation according to \cite{KM2015}. This equation is obtained for a finite box in the physical space while its limiting case for the box size tending to infinity yields the wave kinetic equation.
This approach have been applied for classical NLS (cubic nonlinearity) and to the Hasegawa-Mima equation (quadratic nonlinearity). The main difference between these two cases was that an effective equation in the first case includes all resonances existing in the system while in the latter case it is splitted into a few effective equations (due to the existence of  independent resonance clusters) which can be studied independently,  \cite{KM2016}. Our main goal in this Letter is to deduce the effective equation  for the   Schr\"{o}dinger equation with quadratic nonlinearity (qNLS) and to study the corresponding resonance clustering.

\section{Construction of an effective equation}
We consider the 2D NLS with a quadratic nonlinearity $g$, under the viscosity damping $f(-\Delta)$, where $f\geq 0,$ and a stochastic force in the periodic case with space variables in $\Z_L^2=(\R/2\pi L \Z) \times (\R/2\pi L \Z)$, i.e. the following equation
\be\label{E1}
\partial_t u -i \Delta u = - i \nu \mu g(u, \bar{u}) - \nu f(-\Delta) u + \sqrt{\nu} \frac{d}{dt} \sum_{\ki \in \Z_L^2} b_{\ki} \beta_\ki (t) e^{i \ki \cdot x},
\ee
where $\nu$ is the small parameter, $\mu>0$ is the dispersion coefficient, and $\beta_\ki(t)$ are independent standard complex Wiener processes. Here $L$ is the  size of the box in the $x$-space.

Below we give a brief sketch of the deduction of effective equation  as a consequence of four main steps and final answer for each case of monomial quadratic nonlinearity $g=  u^2, \, \bar{u}u$ and $ \bar{u}^2$.
\begin{itemize}
\item [Step 1.] Recast the PDE into the Fourier modes, and  get a infinite dimensional dynamic systems.

\item [Step 2.] Change the time scale to slow time scale by taking $\tau= \nu t$ as the time variable for the dynamic system. The most important at this step is the demonstration of the well-posedness of the original equation in slow time (which is a nontrivial task), since only in this case, the effective equation is quite close the original equation with small parameter $\nu$, see e.g. \cite{HKM15}, \emph{Theorem 1} for the deterministic case and \emph{Theorem 2} for the stochastic case.

\item [Step 3.] Apply the interaction representation like $a_k= e^{i \nu^{-1} \lambda_k \tau} v_k,$  for including the big linear part (which is of order of $\approx \nu^{-1}$) into the rotation part of the system.

\item [Step 4.] By letting $\nu$ go to zero,  get the effective equation with only the resonant wave vectors left in the equation. This is the so called averaging principle, which means (in a very simplified formulation) that while calculating the limit  $\lim_{T \rightarrow \infty} \frac1T \int _0^T e^{i \theta \tau} d\tau$ (remember now $T \approx \nu^{-1}$) the only nontrivial case is the case where $\theta=0$ and this exactly corresponds to our resonant condition, see e.g. \cite{BM61}, \cite{AKN06} for complete theory.
\end{itemize}

\begin{itemize}

\item[\textbf{Remark}]
The last term $\sqrt{\nu} \frac{d}{dt} \sum_{\ki \in \Z_L^2} b_{\ki} \beta_\ki (t) e^{i \ki \cdot x}$ of \eqref{E1} represents a stochastic force which gives  energy of the system. The parameter $\sqrt{\nu}$ is added  for balancing the friction, since in this case after passing to the slow time at the  Step 2, the stochastic force in the new time scale has the same distribution as the one before.

\end{itemize}

Below we consider all three possible choices of the monomial quadratic nonlinearity and present the result of our calculation:\\
\emph{Case A: nonlinearity $g(u, \bar{u}) =\bar{u}^2$.}\\
Step 1A:
$
\partial_t v_{\ki_3} + i |\ki_3 |^2 v_{\ki_3} = - i \nu \mu \sum_{\ki_1+ \ki_2 + \ki_3=0} \bar{v}_{\ki_1} \bar{v}_{\ki_2} - \nu f(|\ki_3|^2 ) v_{\ki_3} + \sqrt{\nu} b_{\ki_3} \dot{\beta}{\ki_3} (t).
$\\
Step 2A:
$
\partial_t v_{\ki_3} + i \nu^{-1} |\ki_3 |^2 v_{\ki_3} = - i   \mu \sum_{\ki_1+ \ki_2 + \ki_3=0} \bar{v}_{\ki_1} \bar{v}_{\ki_2}  -  f(|\ki_3|^2 ) v_{\ki_3} +   b_{\ki_3} \dot{\beta}_{\ki_3} (\tau).
$ (${\beta}_{\ki_3}$ is another sequence of independent standard complex Wiener processes)\\
Step 3A:
$
\dot{a_{3}}= - i \mu \sum_{\ki_1+ \ki_2 + \ki_3=0 } e^{i \nu^{-1} (|\ki_1|^2 +|\ki_2|^2 +|\ki_3|^2) \tau } \bar{a}_1 \bar{a}_2 - f(|\ki_3|^2) a_3 + b_3 \dot{\beta}_3.
$ (Here and below we write $a_i, b_i, \beta_i$ for $a_{\ki_i}, b_{\ki_i}, \beta_{\ki_i}$, $i=1,2,3,$ for simplicity)\\
Step 4A: In this specific case, we obtain the following degenerate effective equation
\be \label{EE-1}
\dot{a}_3= - f(|\ki_3|^2) a_3 + b_3 \dot{\beta}_3.
\ee
which can be solved  explicitly:
\be
{a}_3(\tau) = a_3(0) e^{-f(|\ki_3|^2)\tau} + b_3 \int_0^\tau e^{-f(|\ki_3|^2) (\tau-r)} d \beta_3 (r).
\ee

\emph{Case B: nonlinearity $g(u, \bar{u}) =u^2$ }\\
Step 3B:
$
\dot{a_{3}} = - i \mu \sum_{\ki_1+ \ki_2 = \ki_3 } a_1 a_2 e^{i \nu^{-1} (-|\ki_1|^2 - |\ki_2|^2 +|\ki_3|^2) \tau }
 - f(|\ki_3|^2) a_3 + b_3 \dot{\beta}_3.
$\\
Step 4B:
\be \label{EE-2}
\dot{a}_3(\tau) = - f(|\ki_3|^2) a_3 - i\mu \sum_{\ki_1, \ki_2 \in S_L^1(\ki_3)} a_1 a_2 + b_3 \dot{\beta}_3,
\ee
where $S_L^1(\ki_3)$ is the resonant set
\be \label{RC1}
S_L^1(\ki_3)=\{(\ki_1, \ki_2) \in \Z_L^2 \times \Z_L^2: \ \ki_1+\ki_2=\ki_3, |\ki_1|^2 +|\ki_2|^2 =|\ki_3|^2\}.
\ee

\emph{Case C: nonlinearity $g(u, \bar{u}) =\bar{u}u$ }\\
Step 3C:
$\dot{a_{3}} = - i \mu \sum_{\ki_1- \ki_2 = \ki_3 } \bar{a}_2 a_1 e^{i \nu^{-1} (|\ki_2|^2 - |\ki_1|^2 +|\ki_3|^2) \tau }
 - f(|\ki_3|^2 ) a_3 + b_3 \dot{\beta}_3.$\\
Step 4C:
\be \label{EE-3}
\dot{a}_3(\tau) = - f(|\ki_3|^2) a_3 - i\mu \sum_{\ki_1, \ki_2 \in S_L^2(\ki_3)} \bar{a}_2 a_1 + b_3 \dot{\beta}_3,
\ee
where $S_L^2(\ki_3)$ is resonant set
\be \label{RC2}S_L^2(\ki_3)=\{(\ki_1, \ki_2)\in \Z_L^2 \times \Z_L^2 : \ \ki_2 +\ki_3=\ki_1, |\ki_2|^2 +|\ki_3|^2 =|\ki_1|^2\}.
\ee
Summing up, we have the following theorem.

\emph{Theorem 1.} \emph{Regard the Eq.~\ref{E1} with three quadratic monomial nonlinearities  $ u^2, \, \bar{u}u, \, \bar{u}^2$. Then corresponding effective equations  read}

\emph{(a) $g=\bar{u}^2:\,\,\, $
${a}_3(\tau) = a_3(0) e^{-f(|\ki_3|^2)\tau} + b_3 \int_0^\tau e^{-f(|\ki_3|^2) (\tau-r)} d \beta_3 (r)$  (no resonances),}

\emph{(b) $g=u^2:\,\,\, $
$\dot{a}_3(\tau) = - f(|\ki_3|^2) a_3 - i\mu \sum_{\ki_1, \ki_2 \in S_L^1(\ki_3)} a_1 a_2 + b_3 \dot{\beta}_3$,}

\emph{(c) $g=\bar{u}u:\,\,\, $
$\dot{a}_3(\tau) = - f(|\ki_3|^2) a_3 - i\mu \sum_{\ki_1, \ki_2 \in S_L^2(\ki_3)} \bar{a}_2 a_1 + b_3 \dot{\beta}_3, $
}

\emph{where two non-empty resonance sets $S_L^1(\ki_3)$ and $S_L^2(\ki_3)$ are given by the Eqs.~(\ref{RC1}),(\ref{RC2}).}

\section{Structure of resonances}
\emph{Resonance set $S_L^1$}.  Resonance conditions
\be \label{RC}
|\ki_1|^2+|\ki_2|^2=|\ki_3|^2, \, \ki_1 + \ki_2= \ki_3,
\ee
imply that $\ki_1 \cdot \ki_2 =0$, i.e.vectors $\ki_1$  and $\ki_2$ are orthogonal. Accordingly, each solution, if any, is formed by three points in $\Integer^2$ and allows a simple geometrical presentation  as a right angle triangle. In order to see whether (\ref{RC}) have a solution, we rewrite it in the coordinates presentation using notations $\ki_j=(k_{j,x}, k_{j,y}), \, j=1,2,$ and deduce that
 \be \label{solutions_general}
k_{1,x} \cdot k_{2,x} = -k_{1,y} \cdot k_{2,y}
 \ee

Resonance cluster is by definition a set of solutions to the resonance conditions having one or more joint wave vectors. Clusters may have different forms, be finite or infinite, have an integrable or non-integrable dynamics, depending on the type of connections between  resonances within a cluster. In a 3-wave systems there are three types of connections: PP (connecting mode is not a high-frequency mode in both solutions), AP (connecting mode is a high-frequency mode in one of the solutions) and AA (connecting mode is a high-frequency mode in both solutions); for details  see \cite{CUP}. In our notations, wavevector $\ki_3$ corresponds to the high-frequency mode. Speaking very generally, three cases are possible: no resonances, all vectors are resonant, and some vectors are resonant and some - not, \cite{PRL}. The following simple properties of the resonance set in our case can be deduced.
\begin{itemize}
\item If $k_{1,x}=0$, it follows from (\ref{solutions_general}) that $k_{2,x}=k_{3,x}$ and $k_{1,y} \cdot k_{2,y}=0$, i.e. whether $k_{1,y} =0$ (degenerate case) or $k_{2,y}=0$. The latter case produces 2-parametric series of solutions $\{(0,k_{1,y}), (k_{2,x},0), (k_{2,x},k_{1,y})\}, \, k_{2,x},k_{1,y} \in \Z$ and each solution
is  presented in the Fourier space by a right angle triangle  in which  catheti lie on the coordinates axes direction. This is a 2-parametric series of solutions and it follows that by fixing one parameter, say $k_{1,y}$, solutions generated by this series form a resonance cluster in which (1) all solutions have one joint vector $(0,k_{1,y})$, and (2) it does not correspond to the high-frequency mode. The clusters of this form containing $N$ triads are called $N$-stars and have an integrable dynamics (with polynomial integrals of motion), under a simple set of assumptions (see \cite{CUP}, p.110, for details). As indexes $1$ and $2$ are interchangeable, the cases when $k_{2,x}$  or $k_{1,y}$  or $k_{2,y}$  are equal to zero are not regarded separately.

\item If $k_{j,x}, k_{j,y} \neq 0 \, \forall j=1,2$, the relation (\ref{solutions_general}) can be regarded as a 3-parametric series of solutions to the resonance conditions (\ref{RC}). Indeed, let
\be \label{3-parametric-series}
k_{2,y} = - k_{1,x} \cdot k_{2,x}/k_{1,y}
\ee
and  choose $k_{1,x}=1, \, k_{2,x}=1, \, k_{1,y}=1$, then the simplest solution to (\ref{RC}) reads
\be \label{solution_simplest}
\ki_1=(1,1), \ki_2=(1,-1), \ki_3=(2,0).
\ee
In particular, each wavector in this case is also a part of an $N$-star cluster. The reason we choose to show the solution (\ref{solution_simplest}) here is that it allows to connect resonance clusters computed from two different series of solutions \emph{via} vector $(2,0)$. Thus we have shown that there are infinitely many solutions to the resonance conditions (\ref{RC}) in the Fourier space and that each wavevector is part of a resonant triad.

    \item All solutions above have been found for the periodic box in space variables with $L$ being  the box's sides, i.e. sides' ratio is 1. The question which is interesting for applications is how the resonant set will be changed if the box sides will be $p$ and $q$ such that $p \neq q, \ p, q, \in \Z^+ $ and for simplicity, they have no common factors, i.e. $(p,q)=1$. For this case it is known that (a) the number of resonances within a periodic box depends drastically on this ratio, and (b) solutions in the periodical box $p, \, q$ are a subset of the solutions in the square box satisfying  suitable divisibility properties. For instance, for all three wave vectors forming a solution,  all $x$-coordinates should be divisible by $p$ and all $y$-coordinates should be divisible by $q$, and the number of boxes can be estimated when  sides ratio $L$ is rational and no resonances occur, \cite{PHD2}. Obviously, if the sides' ratio $L$ is an irrational number, the resonance set is also empty.
\end{itemize}

\emph{Resonance set $S_L^2$.}\\
It is easy to see that the resonance set $S_L^2$ can be transformed into the resonance set $S_L^1$ by switching the indexes $1$ and $3$, i.e. all the properties formulated in the previous Section still hold. The reason why we used these two notations is  as follows. As our next step we plan to study resonances for any linear  combinations of  the monomial quadratic nonlinearities $ u^2, \, \bar{u}u, \, \bar{u}^2$  in which case the position of the running index 3 would be important.

\section{Brief discussion}
We regarded the   Schr\"{o}dinger equation with quadratic nonlinearity  and deduced
effective equation  explicitly for each  case of monomial quadratic nonlinearities $ u^2, \, \bar{u}u, \, \bar{u}^2$. We demonstrated that the set of resonances is empty if the nonlinearity has the  form $\bar{u}^2$ and the effective equation degenerates into one linear stochastic differential equation (\ref{EE-1}) which is solved explicitly. For two other quadratic monomials, the sets of resonances are not empty and corresponding effective equations (\ref{EE-2}), (\ref{EE-3}) are nontrivial.  We have also shown that in these two cases,  \emph{each wave} takes part in an infinite number of resonances, i.e. is a part of an infinite resonance cluster in the Fourier space. Whether there exist independent clusters is still an open problem to be solved in the future as well as specifics of a more general quadratic nonlinearity $u^2 +2 \bar{u}u$ which makes the system Hamiltonian.

The effective equation method allows to regard a more physically relevant setting (by the  including of an external forcing) and thus enrich the classical kinetic WTT. Another important point is that the use of this method allows to deduce (generally infinite) dynamical system for the amplitudes of the resonantly interacting waves while in the discrete WTT dynamical systems are constructed by the way of asymptotical methods as a reduction of an original PDE, for specially chosen resonance condition. It might be that some dynamical characteristics of the wave system are lost. For instance, if one vector takes part in two resonances satisfying different conditions, say, for nonlinearity $ u^2 + 2\bar{u}u$, one solution belongs to $S_L^1$ and another solution belongs to $S_L^2$. Clusters of this kind can not be deduced by the standard asymptotical reduction used in the discrete WTT will be discovered by the effective equation method. Last but not least. The limiting form of the effective equation (with the size of the box going to infinity) is a wave kinetic equation with random data as in  \cite{KM2015,DK19}, or its discrete analog that "does not
involve any randomization of the data", \cite{FGH2016}. Aiming in our future work to deduce a kinetic equation  relevant  for applications, we have chosen to follow the approach developed in  \cite{KM2015,DK19}, with stochastic force introduced according to e.g. \cite{ZL1975,ZLF92}.

\acknowledgments
H. Z. and E. T. acknowledges the support of the Austrian Science Foundation (FWF) under projects P30887 and P31163.
H. Z. is supported by the National Postdoctoral Program for Innovative Talents, No:~BX20180075 and China Postdoctoral Science Foundation.

\end{document}